\documentclass[12pt]{amsart}
\usepackage{amssymb}
\usepackage{amsfonts}
\usepackage{latexsym}

\newcounter{cs}
\stepcounter{cs}
\newcounter{ds}
\stepcounter{ds}

\newcommand{\casos}{\begin{itemize}}
\newcommand{\fcasos}{\end{itemize}\setcounter{cs}{1}}
\newcommand{\cas}{\item[$(\alph{cs})$]\stepcounter{cs}}

\newfont{\got}{eufm10 scaled \magstep1}
\newcommand{\gc}{\mbox{\got c}}
\newcommand{\goe}{\mbox{\got E}}

\newcommand{\ld}{{\rm lim}_{_{\kern-14pt \longrightarrow \kern3pt}}}
\newcommand{\Moo}{M_{\infty}}
\newcommand{\C}{$C^*$-algebra}
\newcommand{\Cs}{$C^*$-algebras}
\newcommand{\del}{\delta}
\newcommand{\eps}{\epsilon}
\newcommand{\alp}{\alpha}
\newcommand{\mul}{{\mathcal M}(A)}
\newcommand{\fl}{\rightarrow}
\newcommand{\sns}{\Leftrightarrow}
\newcommand{\ol}{\overline}
\newcommand{\W}{W_{\sigma}^d(S_u)}
\newcommand{\de}{{{\partial }_{e}}}
\newcommand{\wt}{\widetilde}
\newcommand{\Aff}{\mbox{Aff }}
\newcommand{\LA}{{\rm LAff}}

\newtheorem{lem}{Lemma}[section]
\newtheorem{corol}[lem]{Corollary}
\newtheorem{theor}[lem]{Theorem}
\newtheorem{prop}[lem]{Proposition}
\newtheorem{rema}[lem]{Remark}
\newtheorem{defi}[lem]{Definition}

\begin{document}
\title[Extremal richness of $C^*$-algebras]{Extremal richness of
  multiplier and corona algebras of simple \boldmath $C^*$-algebras
  with real rank zero}
\author{Francesc Perera}
%\thanks{}
\address{Departament de Matem\`atiques, Universitat Aut\`onoma de
  Barcelona, 08193, Bellaterra (Barcelona), Spain}
\email{perera@mat.uab.es}
\date{}
\dedicatory{}
\commby{}

\begin{abstract} 
In this paper we investigate the extremal richness of the multiplier
algebra $\mul$ and the corona algebra $\mul/A$, for a simple \C\ $A$
with real rank zero and stable rank one. We show that the space of
extremal quasitraces and the scale of $A$ contain enough information
to determine whether $\mul/A$ is extremally rich. In detail, if the
scale is finite, then $\mul/A$ is extremally rich. In important cases,
and if the scale is not finite, extremal richness is characterized by
a restrictive condition: the existence of only one infinite extremal
quasitrace which is isolated in a convex sense.
\end{abstract}

\maketitle

\begin{center}
\end{center}

\begin{center}
\end{center}

\section*{Introduction}

The class of \Cs\ with extremal richness was introduced by Brown and
Pedersen in \cite{bpe}, with the objective of extending the theory and
results of finite \Cs\ to the infinite case. Examples of extremally
rich \Cs\ include stable rank one algebras, von Neumann algebras and
purely infinite simple \Cs. Moreover, this class is invariant under
the passage to hereditary subalgebras and under natural constructions
such as tensoring with $M_n(\mathbb{C})$ for all $n\in \mathbb{N}$.

Extremally rich \Cs\ have significant similarities in their properties
with the class of stable rank one \Cs. For example, the presence of
extremal richness gives bounds on the real ranks of the algebras
considered. As is proved in \cite{pelam}, if $A$ is extremally rich,
then $RR(A)\leq 1$. Also, extremally rich \Cs\ with real rank zero are
shown to satisfy the weak cancellation property of separativity (see
\cite{b}, \cite[Theorem 2.11]{bpe4}), that is, the monoid $V(A)$ of Murray-von Neumann
equivalence classes of projections is separative, which means by
definition that whenever $a+a=a+b=b+b$ in $V(A)$ it follows that $a=b$
(see \cite{agop}). Simple and separative \Cs\ with real rank zero are
either purely infinite simple or they have stable rank one
(\cite[Theorem 7.6]{agop}). The same behaviour is observed for simple
\Cs\ with extremal richness (\cite[Corollary 2.8]{bpe2}). This is
therefore closely related to the longstanding Open Question: Is every
finite simple \C\ stably finite?

Our aim in this paper is to analyze extremal richness of multiplier
and corona algebras for a large class of (nonunital) \Cs. We will work
within the class of simple separable \Cs\ with real rank zero and
stable rank one. We also assume that $V(A)$ is strictly
unperforated. This class has been studied in different instances: for
example, see \cite{linpr}, \cite{gok}, \cite{id}; it contains $AF$
algebras, and also many examples which are not $AF$ (\cite{el},
\cite{linpr}, \cite{gpm}). Some of our results will use the additional
hypothesis that the multiplier algebra has real rank zero. As shown by
Lin in \cite[Theorem 10]{linexp}, this occurs if $K_1(A)=0$.

Let $A$ be a nonunital \C\ lying in the abovementioned class, and such
that the real rank of $\mul$ is zero. We will give a complete answer
to the following problem: What conditions determine whether $\mul$ or
$\mul/A$ have extremal richness? For the class of $\sigma$-unital,
purely infinite simple \Cs, for some stabilizations of simple unital
$AF$ algebras, and for simple $AF$ algebras with a finite number of
semi-finite extremal traces, this question has been successfully
considered in \cite{lo}.

Our approach to the problem is based on combining the analysis of
extensions with nonstable $K$-theoretic methods. This involves the
knowledge of properties of $V(\mul)$ and their relation with the ideal
lattice of $\mul$. We will benefit from results concerning this issue,
that appear mainly in \cite{zhriesz}, \cite{linpams}, \cite{gok},
\cite{id}. Thus, the first section is devoted to summarizing the
basics on monoids and their connexion with \Cs\ that will be used in
the sequel. The main objective of section $2$ is to prove the
following: if $A$ has finite scale then $\mul/A$ is extremally rich,
whereas if the scale is infinite, then $\mul/A$ does not have extremal
richness provided that $A$ has at least two infinite extremal
quasitraces. We also discuss basic properties of extremal richness and
some results that will be needed in the last section.

The case in which $A$ has exactly one infinite extremal quasitrace is
handled in section $3$. In this particular situation, the study of
extremal richness of $\mul/A$ requires a deeper analysis, based on the
problem of lifting isometries, for which the index map turns out to be
a useful tool. We close by proving that in our situation $\mul$ is
never extremally rich. It will be clear from this and \cite[Lemma
7.2]{id} that the class of separative \Cs\ properly contains the class
of \Cs\ with extremal richness.

\section{Notation and Preliminaries}

In this section we recall some basic definitions on monoids and \Cs\
that will be used in the subsequent sections. We emphasize the
relation between the order-ideals of the monoid and the closed ideals
of a \C.

All monoids in this paper will be abelian, and consequently we will
write them additively and we will use $0$ for their identity
element. The operation on a monoid $M$ defines a natural preordering
by:
$$x\leq y \sns y=x+z \mbox{ for some } z\in M,$$
which is translation-invariant. This preordering is sometimes called
the {\bf algebraic preordering}. As usual, we write $x<y$ if $x\leq y$
and $x\ne y$.

If $M$ is a  monoid, a nonzero element $u\in M$ is called an {\bf
  order-unit} if for any $x\in M$, there exists $n\in \mathbb{N}$ such
that $x\leq nu$. We say that $M$ is {\bf conical} provided that the
set $M^*$ of nonzero elements is closed under addition. For a \C\ $A$,
we denote by $V(A)$ the monoid of Murray-von Neumann equivalence
classes of projections from $\Moo (A)$. (Equivalently, if $A$ is
unital $V(A)$ can be described as the additive monoid of isomorphism
classes of finitely generated projective modules over $A$.) Note that
$V(A)$ is always conical, and that $[1_A]$ is an order-unit for $V(A)$
if $A$ is unital.

A nonempty subset of a monoid which is a submonoid and
order-hereditary will be called an {\bf order-ideal}. We say that a
monoid $M$ is {\bf simple} if $M$ has precisely two order-ideals,
namely the ideal generated by $0$ and $M$. In case $M$ is conical,
then $M$ is simple if and only if $M$ is nonzero and every nonzero
element is an order-unit. This is the case for $V(A)$, where $A$ is a
simple \C.

Let $M$ be a monoid and let $I$ be an order-ideal of $M$. Define an
equivalence relation in $M$ as follows: if $x,y\in M$ write $x\sim y$
if and only if there exist $z,w\in I$ such that $x+z=y+w$. Denote by
$M/I$ the quotient of $M$ modulo this equivalence relation, and by
$[x]$ the equivalence class of an element $x\in M$. The addition
$[x]+[y]=[x+y]$ is then a well-defined operation under which $M/I$
becomes a monoid, referred as to the {\bf quotient monoid} of $M$
modulo $I$. If $A$ is a \C\ and $I$ is a closed ideal of $A$, then
$V(I)$ is naturally an order-ideal of $V(A)$. Moreover, if $A$ has
real rank zero, then the quotient $V(A)/V(I)$ is isomorphic to
$V(A/I)$ (\cite[Proposition 1.4]{agop}).

We say that a cancellative monoid $M$ is {\bf strictly unperforated}
if whenever $nx<ny$ for some $n\in \mathbb{N}$ and $x,y\in M$, it
follows that $x<y$. It is remarkable that no examples are known of
simple \Cs\ $A$ with real rank zero and stable rank one whose $V(A)$'s
are not strictly unperforated, and therefore this technical condition
is quite natural. (If the real rank zero condition is dropped, then
there are examples with perforation on $V(A)$, as shown in \cite{vd}.)

Let $K$ be a compact convex set. We denote by ${\rm LAff}(K)$ the
monoid of all affine and lower semicontinuous functions on $K$ with
values on $\mathbb{R}\cup\{\infty\}$, and we shall use ${\rm Aff}(K)$
to denote the submonoid of elements in ${\rm LAff}(K)$ that are
continuous. Let ${\rm LAff}_{\sigma}(K)$ be the submonoid of ${\rm
  LAff}(K)$ whose elements are pointwise suprema of increasing
sequences of elements from ${\rm Aff}(K)$. The use of the superscript
$++$ will always refer to strictly positive functions.

\section{Extremal richness of multiplier and corona algebras}

In this section we introduce the class of extremally rich \Cs, giving
some equivalent definitions, and discussing related matters about
extensions that can be found in \cite{bpe} and \cite{lo}. We present
at the end a first result that analyzes the extremal richness of the
multiplier and corona algebras, for a wide class of simple \Cs\ with
real rank zero and stable rank one. To establish this fact, we need
some results concerning the ideal structure of these rings, that
appear in \cite{id}, and thus they will be stated as required.

If $A$ is a unital \C, we use $\goe (A)$ to denote the set of extreme
points of its closed unit ball $A_1$, and we refer to this set as the
set of extreme points of the algebra. Recall that the elements of
$\goe (A)$ are precisely those partial isometries $v\in A$ satisfying
$(1-vv^*)A(1-v^*v)=0$ (see, for example, \cite[Proposition
1.4.7]{ped}). Notice that if $A$ is prime, then the extreme points are
precisely the isometries and co-isometries of the algebra. An element
$x\in A$ is said to be {\bf quasi-invertible} if $x\in A^{-1}\goe
(A)A^{-1}$, and the set of quasi-invertible elements is denoted by
$A_q^{-1}$.
\begin{defi}
\cite[Section 3]{bpe}
We say that a (unital) \C\ $A$ is {\bf extremally rich} if the set
$A_q^{-1}$ of quasi-invertible elements is dense in $A$. As usual, a
nonunital \C\ $A$ is extremally rich if its minimal unitization $\wt
{A}$ is extremally rich.
\end{defi}

An equivalent notion may be found in \cite[Section 3]{bpe} (see also
\cite{bpe3}): A unital \C\ $A$ is extremally rich if and
only if $A_1={\rm conv}(\goe (A))$. At this point, it is convenient to
notice that in any unital \C\ $A$, the {\it closure} of the convex
hull of the unitaries of $A$ (which are extreme points) equals $A_1$
(see, e.g. \cite[Proposition 1.1.12]{ped}). We denote by $\mathcal{U}
(A)$ the group of unitaries of a (unital) \C\ $A$.
\begin{rema}
\label{rm}
Let $A$ be a \C. Then:
\casos
\cas
{\rm (\cite[Section 1]{b}, \cite[Proposition 2.6]{bpe2})} $A$ has stable rank one if and only if $A$
is extremally rich and $\goe (\wt{A})=\mathcal{U}(\wt{A})$.
\cas
{\rm (\cite{b}, \cite[Corollary 2.8]{bpe2}, \cite[Lemma 3.3]{lo})} If
$A$ is simple, then $A$ is extremally rich if and only if $A$ is
either purely infinite or it has stable rank one.
\fcasos
\end{rema}

A consequence of Remark~\ref{rm} and \cite[Theorem 7.6]{agop} is that
if $A$ is a simple \C\ with real rank zero, then $A$ has extremal
richness if and only if it is separative. This contrasts with the fact
that, as we will see, the class of \Cs\ with extremal richness is
strictly contained in the class of separative \Cs.

As in the case of \Cs\ with real rank zero (see \cite[Theorem
3.14]{bpf}, \cite[3.2]{zhi}), the behaviour of extremal richness under
extensions depends not only on the extremal richness of the ideal and
the quotient algebra, but also on a lifting condition of the extreme
points. We now record some results in this direction:
\begin{theor}
\label{1.1}
\cite[Theorem 6.1]{bpe}
Let $J$ be a closed ideal in a unital \C\ $A$. Then $A$ is extremally
rich if and only if $J$ and $A/J$ are extremally rich, the extreme
points of $A/J$ lift to those of $A$ and $\goe (A)+J\subset
(A_q^{-1})^{-}$.\qed
\end{theor}
\begin{corol}
\label{1.1a}
\cite[Corollary 6.3]{bpe}
Let $A$ be a \C\ and let $J$ be a closed ideal of $A$ with stable rank
one. Then $A$ is extremally rich if and only if $A/J$ is extremally
rich and the extreme points of $A/J$ lift to those of $A$.\qed
\end{corol}
\begin{theor}
\label{1.3a}
\cite[Theorem 3.6]{lo}
Let $A$ be a \C, and let $J$ be an essential closed ideal of $A$ which
is purely infinite simple. Then $A$ is extremally rich if and only if
$A/J$ is extremally rich and $\goe (A/J)$ consists of isometries and
co-isometries.\qed
\end{theor}

Let $M$ and $N$ be monoids with respective order-units $u$ and $v$. A
monoid morphism (that is, an additive map)  $f:M\fl N$ is said to be
{\bf normalized} provided that $f(u)=v$. Recall that a {\bf state} on
a monoid $M$ with order-unit $u\in M$ is a normalized monoid morphism
$s:(M,u)\fl (\mathbb{R}^+,1)$. We denote the set of states on $(M,u)$
by $St(M,u)$ or by $S_u$ when no confusion may arise. We also denote
by $\phi_u:M\fl {\rm Aff}(S_u)$ the natural map, given by evaluation
on the states of $M$. Observe that $St(M,u)=St(G(M),u)$, where $G(M)$
is the Grothendieck group of $M$, and hence it is a compact convex
set.

In order to analyze the extremal richness of the multiplier and corona
algebra of a simple \C\ $A$ with real rank zero and stable rank one,
the ideal structure of $\mul$ will play a crucial role. Let $A$ be a
$\sigma$-unital simple \C\ with real rank zero and stable rank
one. Fix $u\in V(A)^*$, and set $d=\sup\phi_u(D(A))$. We define:
$$\W=\{f\in {\rm LAff}(S_u)^{++}\mid f+g=nd \mbox{ for some } g\in
{\rm LAff}(S_u)^{++} \mbox{ and }n\in \mathbb{N}\}.$$
Consider the set $V(A)\sqcup \W$, where $\sqcup$ stands for disjoint
union of sets. We equip this set with a monoid structure that extends
the natural given addition operations of both $V(A)$ and $\W$, and by
setting $x+f=\phi_u(x)+f$, for $x\in V(A)$ and $f\in \W$. It is not
difficult to see that this is a well defined operation. In
\cite[Theorem 3.10]{id} the following important relation between
$V(\mul)$ and $V(A)$ was established:
\begin{theor}
\label{tm}
Let $A$ be a $\sigma$-unital nonunital \C. Suppose that $A$ is simple,
with real rank zero, stable rank one and that $V(A)$ is strictly
unperforated. Assume that $A$ is nonelementary. Fix a nonzero element
$u\in V(A)$. Set $D(A)=\{[p]\in V(A)\mid p \mbox{ is a projection in
  }A\}$, and $d=\sup\phi_u(D(A))$. Then there is a normalized monoid
isomorphism
$$\varphi: V(\mul)\fl V(A)\sqcup \W,$$
such that $\varphi([p])=[p]$ if $p\in A$, and
$\varphi([p])=\sup\{\phi_u ([q])\mid [q]\in V(A)\mbox{ and } q\lesssim
p\}$ if $p\in \mul\setminus A$.\qed
\end{theor}

We will make an extensive use of this result in the next section; for
now we content ourselves with two applications. For a compact convex
set $K$, denote by $\de K$ the set of its extreme points. Combined
with \cite[Theorem 2.3]{zhriesz} (see also \cite[Theorem 2.1]{id}),
Theorem~\ref{tm} gives an effective method to study the ideal
structure of multiplier and corona algebras of \Cs\ in the class we
are considering. In detail, the map $I\mapsto \varphi(V(I))$ provides
a lattice isomorphism between the lattice of closed ideals of $\mul$
and the lattice of order-ideals of $V(A)\sqcup\W$. In this context, we
define the {\bf finite ideal} of $\mul$ as the unique closed ideal
$I_{fin}(A)$ of $\mul$ such that $\varphi (V(I_{fin}(A)))=V(A)\sqcup
\{f\in \W\mid f_{|\de S_u} \mbox{ is finite}\}$.

Let $A$ be a simple \C\ with real rank zero, and let $u\in V(A)^*$. We
say that $A$ has {\bf finite scale} provided that the (lower
semicontinuous) affine function $d:=\sup\phi_u(D(A))$ is finite when
restricted to $\de S_u$ (see \cite[Definition 4.6]{id}). It should be
noted that this notion does not depend on the choice of the nonzero
element $u\in V(A)^*$, and that it differs from the definition of
finite scale given in \cite{linaf}, in that the condition on $d$ is
required only for the extreme boundary of the state space. In our
setting, \Cs\ with finite scale are characterized by the following
nice property:
\begin{theor}
\label{tm2}
{\rm (cf. \cite[Theorem 4.8]{id})}
Let $A$ be a nonunital simple and separable \C\ with real rank zero
and stable rank one. Assume that $A$ is nonelementary, that $V(A)$ is
strictly unperforated and that $\mul$ has real rank zero. Then $A$ has
finite scale if and only if, for every closed ideal $I$ of $\mul$
properly containing $A$, we have that ${\rm sr}(\mul/I)=1$.\qed
\end{theor}
\begin{defi}
\label{qt}
\cite[Definition 5.1]{id}
Let $A$ be a \C. A {\bf\boldmath $1$-quasitrace} on $A$ is a map
$\tau:A_+\fl [0,\infty]$ such that $\tau(\alp x)=\alp\tau(x)$ if $x\in
A_+$ and $\alp\in \mathbb{R}^+$, such that $\tau
(x+y)=\tau(x)+\tau(y)$, whenever $x$ and $y$ are commuting elements in
$A_+$, and such that $\tau(xx^*)=\tau(x^*x)$ for all $x\in A$. A {\bf
  quasitrace} on $A$ is a $1$-quasitrace $\tau$ that extends to a
$1$-quasitrace $\tau_n$ on $M_n(A)$ for each $n\in \mathbb{N}$.
\end{defi}

We use the convention here that $0\cdot \infty=0$ so that
$\tau(0)=0$. Viewing $A$ as the upper left hand corner subalgebra of
$M_n(A)$, the extension $\tau_n$ in \ref{qt} of $\tau$ means that
$\tau(x)=\tau_n(xe_{11})$, where $e_{11}$ is the matrix unit in
$M_n(\wt{A})$. If $\tau$ is a quasitrace, we say that $\tau$ is {\bf
  densely defined} provided that the set $F_{\tau}:=\{x\in A_+\mid
\tau(x)<\infty\}$ is dense in $A_+$. We denote the set of densely
defined quasitraces by $QT_d(A)$, and we also use $LQT(A)$ to denote
the set of lower semicontinuous quasitraces. The notation $LQT_d(A)$
will stand for the set of lower semicontinuous, densely defined,
quasitraces. If $x\in K(A)_+$, where $K(A)$ is the Pedersen ideal of
$A$, we set $Q=\{\tau\in LQT_d(A)\mid \tau_{|K(A)}<\infty\}$ and
$Q_x=\{\tau\in Q\mid \tau(x)=1\}$. If $A$ is simple, the set $Q_x$ is
(weakly) compact. If $A$ is moreover $\sigma$-unital with real rank
zero and $u=[p]\in V(A)^*$ for a nonzero projection $p\in A$, then the
natural map $\alp:Q_p\fl S_u$ given by $\alp(\tau)([q])=\tau(q)$, for
$[q]\in V(A)$, provides an affine homeomorphism (\cite[Theorem
5.6]{id}). This is the $\sigma$-unital, nonunital and semi-finite
version of Blackadar and Handelman's theorem (\cite[Theorem
III.1.3]{bh}).
\begin{defi}
Let $A$ be a \C. A lower semicontinuous and order-pre\-ser\-ving
quasitrace $\tau$ is said to be {\bf infinite} if
$\sup\limits_{\lambda}\tau(u_{\lambda})=+\infty$ for some approximate
unit $\{u_{\lambda}\}_{\lambda\in \Lambda}$.
\end{defi}

Observe that this definition does not depend on the particular
approximate unit. If $A$ has real rank zero, and $\tau$ is infinite,
then $\sup\tau(p)=\infty$, where the supremum is taken over the
projections $p\in A$.
\begin{theor}
\label{tm3}
{\rm (cf. \cite[Theorem 6.3, Proposition 6.5, Theorem 6.6]{id})}
Let $A$ be a separable nonunital simple \C\ with real rank zero,
stable rank one and with $V(A)$ strictly unperforated. Suppose that
$A$ is nonelementary. Let $p\in A$ be a nonzero projection and let
$\gc$ be the cardinal of infinite extremal quasitraces in $Q_p$. Then
$\mul$ has at least $\gc$ different maximal ideals that contain
$I_{fin}(A)$, and the quotient of $\mul$ by any of these ideals is a
purely infinite simple \C.\qed
\end{theor}

In order to clarify the exposition, we state a lemma whose argument is
used in \cite[Theorem 4.9]{lo}.
\begin{lem}
\label{1.2}
Let $A$ be a unital \C. Assume that $A$ is prime and that there exist
different maximal ideals $J_i$, for $i=1,2$ such that each projection
in $A/J_i$ is infinite. Then $A$ is not extremally rich.
\end{lem}
{\bf Proof.}
Consider the $C^*$-exact sequence:
$$0\fl (J_1\cap J_2)\fl A\fl A/(J_1\cap J_2)\fl 0,$$
and note that $A/(J_1\cap J_2)\cong J_1/(J_1\cap J_2)\oplus
J_2/(J_1\cap J_2)$ and that $J_1/(J_1\cap J_2)\cong A/J_2$ and
$J_2/(J_1\cap J_2)\cong A/J_1$. Since each quotient $A/J_i$ is simple,
its extreme points are isometries or co-isometries. Further, as every
projection in $J_i/(J_1\cap J_2)$ is infinite, we see that
$J_i/(J_1\cap J_2)$ contain non-trivial isometries (and hence
co-isometries) for $i=1,2$. Finally, since the set of extreme points
of a direct sum equals the direct sum of the extreme points of each
factor, we see that $A/(J_1\cap J_2)$ has an extreme point which is
neither an isometry nor a co-isometry, and hence it cannot be lifted
to an extreme point of $A$, because $A$ is prime. By Theorem
\ref{1.1}, $A$ is not extremally rich.\qed

Let us now discuss the extremal richness in the elementary
case. Suppose that $A\cong \mathbb{K}$, where
$\mathbb{K}=\mathbb{K}(\mathcal{H})$ is the $C^*$-algebra of compact
operators over an infinite-dimensional, separable, Hilbert space
$\mathcal{H}$. Then, if $\mathbb{B}=\mathbb{B}(\mathcal{H})$, we have
that $\mul/A\cong \mathbb{B}/\mathbb{K}$ is a purely infinite simple
\C, hence extremally rich by Remark~\ref{rm} $(a)$. On the other hand,
$\mathbb{B}$ is a von Neumann algebra, whence it is also extremally
rich. Thus we shall assume from now on that all \Cs\ are
nonelementary.

In \cite[Theorem 4.9]{lo} it is proved that if $A$ is a simple
separable $AF$ algebra such that $A\otimes \mathbb{K}$ contains at
least two semi-finite extremal traces, then $\mathcal{M}(A\otimes
\mathbb{K})/(A\otimes \mathbb{K})$ is not extremally rich. On the
other hand, in \cite[Proposition 4.13]{lo} it is established that if
$A$ is a simple and separable (nonunital) $AF$ algebra with a finite
number of semi-finite extremal traces, of which at least two are
infinite, then $\mul/A$ is not extremally rich. Both situations can be
handled in our setting. The next result enlarges to a great extent the
number of instances in which the extremal richness of the corona
algebra can be analyzed. We also answer an implicit question that is
posed in \cite[Remark 4.19]{lo}: If $A$ has an infinite number of
extremal quasitraces, has the corona extremal richness? In case none
of them are infinite, the answer is positive (see also \cite[Theorem
4.1]{lo}, where some stable cases outside our class are considered),
whereas if at least two of them are infinite, the answer is
negative. The case in which there is only one infinite extremal
quasitrace will be considered in the next section.
\begin{theor}
\label{1.4}
Let $A$ be a nonunital separable simple \C\ with real rank zero and
stable rank one. Suppose that $A$ is nonelementary and that $V(A)$ is
strictly unperforated. Let $p\in A$ be a nonzero projection.
\casos
\cas
If $A$ has finite scale and $\mul$ has real rank zero, then $\mul /A$
is extremally rich.
\cas
If $A$ has at least two infinite extremal quasitraces in $Q_p$, then
$\mul /A$ is not extremally rich. In particular, $\mul$ is not
extremally rich.
\fcasos
\end{theor}
{\bf Proof.}
$(a)$. Recall that $\mul$ has a unique closed ideal $L(A)$ that
properly contains $A$ and that is contained in every closed ideal that
properly contains $A$ (see \cite[Remark 2.9]{linpams}, and also
\cite[Proposition 4.1]{id}). Notice that $L(A)/A$ is an essential
closed ideal of $\mul/A$, and that it is purely infinite simple (see
\cite[Theorem 1.3 (a)]{zhcan}). Since $A$ has finite scale, we have
that $\mul /L(A)$ has stable rank one, by Theorem~\ref{tm2}. Using
Remark~\ref{rm} $(a)$, we conclude that $\mul/L(A)$ is extremally rich
and that $\goe (\mul/L(A))=\mathcal{U} (\mul/L(A))$. By
Theorem~\ref{1.3a}, it follows that $\mul/A$ is extremally rich.

$(b)$. Let $\gc$ be the cardinal of infinite extremal quasitraces in
$Q_p$. By hypothesis $\gc\geq 2$. Using Theorem~\ref{tm3}, we get at
least $\gc$ different closed maximal ideals in $\mul/A$. Moreover, the
quotient of $\mul/A$ by each one of these ideals is a purely infinite
simple \C. Also, since $L(A)/A$ is the minimal nonzero closed ideal of
$\mul/A$, we get that the corona algebra is a prime ring. Therefore,
the hypotheses of Lemma~\ref{1.2} are fulfilled, whence we conclude
that $\mul/A$ is not extremally rich.

That $\mul$ is not extremally rich follows from Theorem~\ref{1.1}.\qed

\begin{rema}
{\rm In the proof of the previous theorem, we used the fact that the
  corona algebra $\mul/A$ is a prime ring. Although it is possibly
  well-known, we remark that this is not true in general. Let $A$ be a
  simple \C\ such that $\mul $ has at least two different maximal
  closed ideals $I_1$, $I_2$ (such examples exist; see \cite{id}). Let
  $J=I_1I_2$. Then $J$ is a prime \C\ and $\mathcal{M}(J)=\mul
  $. Therefore ${\mathcal M}(J)/J$ contains two nontrivial ideals,
  $I_i/J$, $i=1,2$, whose product is zero.}
\end{rema}

\section{Extremal richness of corona algebras with only one infinite
  extremal quasitrace}

The purpose of the present section is to determine when the corona
algebra of a simple \C\ $A$ with real rank zero, stable rank one, with
$V(A)$ strictly unperforated, and that has precisely one infinite
extremal quasitrace is extremally rich. Our approach to the solution
follows the lines indicated by Theorem~\ref{1.3a} and
Corollary~\ref{1.1a}. The knowledge of significant aspects of the
ideal lattice of $\mul/A$, of projections in $\mul$ and in its
quotient algebras will be an important ingredient in the
following. This will be reflected in some computations of the index
map in various situations, for which Theorem~\ref{tm} will be
essential.
\begin{lem}
\label{2.1}
Let $B$ be a (unital) \C, and let $I$ be a closed two-sided ideal of
$B$. Let $w$ be an isometry in $B/I$, and denote by $\pi:B\fl B/I$ the
natural quotient map. Then $w$ can be lifted to an isometry $z\in B$
if and only if there exists a partial isometry $v\in B$ such that
$1-v^*v\lesssim 1-vv^*$ and $\pi (v)=w$.
\end{lem}
{\bf Proof.}
If $w$ can be lifted to an isometry $z\in B$, then just take
$v=z$. For the converse, assume that $w=\pi(v)$ for some partial
isometry $v\in B$ such that $1-v^*v\lesssim 1-vv^*$. Let $p=v^*v$ and
$q=vv^*$. Then there exists $r\in B$ such that $1-p=r^*r$ and
$rr^*\leq 1-q$. Since $\pi (v)=w$, we see that $\pi (r^*r)=0$, that
is, $r\in I$. Let $z=v+r$. Then $z^*z=1$, hence $z$ is an isometry,
and $\pi (z)=\pi (v)=w$.\qed

Recall that a monoid $M$ is a {\bf refinement monoid} if whenever
$x_1,x_2,y_1,y_2 \in M$ satisfy $x_1+x_2=y_1+y_2$, then there exist
elements $z_{ij}\in M$, for $i,j=1,2$ such that $\sum\limits_{j}
z_{ij}=x_i$ and $\sum\limits_{i}z_{ij}=x_j$ for each $i,j$. By
\cite[Lemma 2.3]{ap} (see also \cite[Theorem 1.1]{zhriesz}), if $B$ is
a \C\ with real rank zero, then $V(B)$ is a refinement monoid.
\begin{lem}
{\rm (P. Ara)}
\label{g}
Let $B$ be a \C\ with real rank zero. Then $K_0(B)=G(V(B))$.
\end{lem}
{\bf Proof.}
We may clearly assume that $B$ is nonunital. First, notice that
$V(B)=\ld V(pBp)$, where $p$ runs over the set of projections of
$B$. (In general, the set of projections of $B$ need not be directed,
but because of the real rank zero condition it is easy to see that
given projections $p,q\in B$, there exists a projection $r\in B$ such
that $V(pBp), V(qBq)\subseteq V(rBr)$.) Taking into account that $G$
is in fact a continuous functor from the category of monoids to the
category of groups, we get that $G(V(B))=\ld G(V(pBp))$. Since $pBp$
does have a unit for each projection $p\in B$, it follows that
$G(V(pBp))=K_0(pBp)$. Hence, it remains to prove that $K_0(B)=\ld
K_0(pBp)$.

We denote by $B^+:=B\oplus \mathbb{C}$, equipped with pointwise sum
and adjoint, and with a mixed multiplication given by
$(x,\lambda)(y,\mu)=(xy+\mu x+\lambda y,\lambda\mu)$, for $x,y\in B$
and $\lambda,\mu\in \mathbb{C}$. Then $B^+$ is a \C, which is
isomorphic to $\wt{B}$ if $B$ is nonunital (\cite[Proposition
2.1.7]{wo}). Note that if $p\in M_n(B^+)$ is a projection, then there
exist $r\leq n$ and projections $g,h\in \Moo (B)$ such that $g\leq
1_r$ and $p\sim (1_r-g)\oplus h$. (Here $1_r$ stands for the unit of
$M_r(\mathbb{C})$.) This follows from \cite[Lemma 10.3]{gok} (see also
\cite[Lemma 3.4]{pext}). Applying the refinement property to the
equality $(1_r-g)+g=1_r$ we get projections $p_1,\ldots,p_r\in B$ such 
that $g\sim \sum\limits_{i=1}^r \oplus p_i$, while $1_r-g\sim
\sum\limits_{i=1}^r \oplus (1-p_i)$. We therefore obtain that
$1_r-g\sim(1-p_1)\oplus\ldots\oplus (1-p_r)\sim 1_r-\oplus_{i=1}^r
p_i$. Using this and the fact that $V(B)=\ld V(pBp)$, we conclude that
$V(B^+)=\ld V((pBp)^+)$. On the other hand, if $\pi:B^+\fl \mathbb{C}$
and $\pi_p:(pBp)^+\fl \mathbb{C}$, for $p\in B$ are the natural
projection maps, then (using the functoriality of $V$) we have
$V(\pi_p)=V(\pi)\circ V(i)$, where $i:(pBp)^+\fl B^+$ is the natural
inclusion.

Using again the continuity of $G$, we obtain that $K_0(B^+)=\ld
K_0((pBp)^+)$, and also that $G(\pi_p)=G(\pi)\circ G(i)$. We conclude
then that $K_0(B)=\ld K_0(pBp)$, as desired.\qed

As is well known, the previous result is false in general. For
example, let $A=C_0(\mathbb{R}^2)$. Then $V(A)=0$, whereas
$K_0(A)\cong \mathbb{Z}$ (see \cite[6.2]{wo}).

If $B$ is a \C\ and $I$ is a closed ideal of $B$, we denote by $\del
:K_1(B/I)\fl K_0(I)$ the index map in $K$-Theory.
\begin{prop}
\label{2.2}
Let $B$ be a (unital) \C, and let $I$ be a closed ideal of $B$. Let
$w\in \mathcal{U} (B/I)$, and let $\pi:B\fl B/I$ be the natural
quotient map. If $RR(I)=0$ and $V(I)$ is cancellative, then $w$ can be
lifted to a unitary in $B$ (resp. a proper isometry, a proper
co-isometry) if and only if $\del [w]=0$ (resp. $\del [w]<0$, $\del
[w]>0$).
\end{prop}
{\bf Proof.}
If $w$ can be lifted to a partial isometry $v\in B$ via $\pi$, then it
is a standard fact that the index can be computed as $\del
([w])=[1-v^*v]-[1-vv^*]$ (see, for example, \cite[Exercise 8C]{wo}).

Assume that $w$ can be lifted to a unitary (resp. a proper isometry, a
proper co-isometry) $v\in B$. Then it follows easily from the previous
observation (using Lemma~\ref{g} and that $V(I)$ is cancellative) that
$\del [w]=0$ (resp. $\del[w]<0$, $\del[w]>0$).

Conversely, since $RR(I)=0$ there exists a partial isometry $v\in B$
such that $v$ is a lift for $w$ (see the proof of \cite[Lemma
2.6]{elan}, or \cite[Lemma 2.1]{pext}). If $\del[w]=0$, then
$[1-v^*v]=[1-vv^*]$ in $K_0(I)$. By Lemma~\ref{g}, and since $V(I)$ is
cancellative, we get that $1-v^*v\sim 1-vv^*$ in $I$. Therefore, there
exists $r\in I$ such that $1-p=r^*r$ and $1-q=rr^*$, where $p=v^*v$
and $q=vv^*$. Let $z=v+r$. Then $z^*z=zz^*=1$ and $\pi (z)=\pi
(v)=w$. Hence $w$ can be lifted to a unitary.

If $\del[w]<0$, then $1-v^*v\lesssim 1-vv^*$, so that by Lemma
\ref{2.1}, $w$ can be lifted to an isometry $z\in B$, which is proper
since $\del ([w])\ne 0$. We proceed similarly if $\del [w]>0$, in
order to get a co-isometry.\qed

In the following lemma we need to establish a slight generalization of
a known result (\cite[Lemma II.7.1]{al}).
\begin{lem}
\label{pacha}
Let $K$ be a compact convex set, and let $f,g\in {\rm
  LAff}_{\sigma}(K)^{++}$ be functions such that $f_{|\de K}=g_{|\de
  K}$. Then $f=g$.
\end{lem}
{\bf Proof.}
Since $K$ is compact and $f$ is lower semicontinuous, we see that $f$
takes its minimum value on $K$. Indeed, if $\alp=\min (f)$, then there
exists $x\in \de K$ such that $f(x)=\alp$ (this follows after a
standard argument - see \cite[Corollary 5.19]{poag}).

Write $g=\sup_n g_n$, where $g_n\in {\rm Aff}(K)^{++}$ for each $n$,
and observe that $\{g_n\}$ form an increasing sequence. Note that
$f-g_n=g-g_n\geq 0$ on $\de K$ for every $n\in \mathbb{N}$. Taking
into account that $f-g_n$ is affine and lower semicontinuous, we
conclude from the previous paragraph that $f-g_n\geq 0$
globally. Hence $f\geq g$. A similar argument shows that $f\leq g$,
and therefore $f=g$.\qed

The following fact is a consequence of the methods developed in
\cite[Section 4]{id}. Note that if $K$ is a metrizable Choquet
simplex, then ${\rm LAff}_{\sigma}(K)={\rm LAff}(K)$ (this is probably
standard; for a proof, see \cite[Lemma 4.12]{id}).
\begin{lem}
\label{2.3}
Let $K$ be a metrizable Choquet simplex. Let $f,g,h,d\in \LA (K)^{++}$
such that $f+g=f+h=d$. Suppose that there exists $s\in \de K$ such
that $d(x)=\infty$ (for $x\in \de K$) if and only if $x=s$. Assume
that $g_{| \de K}$ is a finite function.
\casos
\cas
If $h_{|\de K}$ is also finite and $g(s)<h(s)$, then there exist $e\in
\LA (K)^{++}$ and $k\in \mathbb{N}$ such that $g+k=h+e$ and
$e+d=k+d$.
\cas
If $h(s)=\infty$ and $\{ s\}'$, the complementary face of $\{s\}$, is
a closed face, then there exist $e\in \LA (K)^{++}$ and $k\in
\mathbb{N}$ such that $g+e=h+k$ and $e+d=k+d$.
\fcasos
\end{lem}
{\bf Proof.}
Since $g_{|\de K}$ is finite, we immediately get that $f(s)=\infty$
and that $g_{|\de K\setminus\{s\}}=h_{|\de K\setminus\{s\}}$.

$(a)$. Let $a=h(s)-g(s)$ and take $k\in \mathbb{N}$ such that
$k>a$. By an argument similar to the one used in \cite[Proposition
4.13]{id}, there exists a lower semicontinuous affine function $e$
such that $e(s)=k-a$ and $e_{|\{s\}'}=k$. Note that $(g+k)_{|\de
  K}=(h+e)_{|\de K}$ and that $(e+d)_{|\de K}=(k+d)_{|\de
  K}$. Therefore $g+k=h+e$ and $e+d=k+d$, by Lemma~\ref{pacha}.

$(b)$. Let $k\in \mathbb{N}$. Since $\{s\}'$ is closed, there is (by
\cite[Corollary 11.27]{poag}) an isomorphism 
$${\rm Aff}(K)\cong {\rm Aff}(\{s\})\times {\rm Aff}(\{s\}').$$
Thus, for each $n\in \mathbb{N}$, there exists $e_n\in {\rm
  Aff}(K)^{++}$ such that $e_n (s)=n$ and ${e_n}_{|\{s\}'}=k$. Let
$e=\sup\limits_n e_n$. Then $e\in {\rm LAff}(K)^{++}$, and satisfies
$e(s)=\infty$ while $e_{|\{ s\}'}=k$. Therefore $(g+e)_{|\de
  K}=(h+k)_{|\de K}$ and $(e+d)_{|\de K}=(k+d)_{|\de K}$, whence
$g+e=h+k$ and $e+d=k+d$, also by Lemma~\ref{pacha}.\qed

\begin{prop}
\label{2.4}
Let $A$ be a (nonunital) simple separable \C\ with real rank zero and
stable rank one. Assume that $A$ is nonelementary and that $V(A)$ is
strictly unperforated. Let $p\in A$ be a nonzero projection and
suppose that $A$ has exactly one infinite extremal quasitrace in
$Q_p$. Suppose also that the real rank of $\mul$ is zero. Let $I$ be
any closed ideal of $\mul$ such that $A\subsetneqq I\subseteq
I_{fin}(A)$. Then $\del ([w])$ is either zero, or positive, or
negative, for any unitary $w\in \mul/I$, where $\del :K_1(\mul/I)\fl
K_0(I/L(A))$ is the index map. Hence $w$ can always be lifted to an
isometry or a co-isometry.
\end{prop}
{\bf Proof.}
Let $u=[p]\in V(A)$, and let $d=\sup\phi_u(D(A))$. By
Theorem~\ref{tm}, there exists a normalized monoid isomorphism
$\varphi: V(\mul )\fl V(A)\sqcup \W$ (so that
$\varphi([1_{\mul}])=d$). Since the natural map $\alp:Q_p\fl S_u$
given by evaluation is an isomorphism, the fact that $A$ has exactly
one infinite extremal quasitrace means that the set $\Gamma_d:=\{t\in
\de S_u\mid d(t)=\infty\}$ consists of one point, which we denote by
$s$.

Let $\phi :\mul \fl \mul /I$ and $\pi:\mul \fl \mul /L(A)$ be the
natural projection maps. By the proof of \cite[Lemma 2.6]{elan} (or
also \cite[Lemma 2.1]{pext}), there exists a partial isometry $v\in
\mul$ such that $\phi (v)=w$. Then $1-v^*v$ and $1-vv^*$ belong to
$I$. Since $w$ is not zero, we first note that $v\notin A$. Suppose
that $1-v^*v\in A$. Then $\del [w]=-[\pi(1-vv^*)]\leq 0$, and a
similar conclusion would result if $1-v^*v\in A$. Therefore we may
assume that $1-v^*v,1-vv^*\notin A$. Let $f=\varphi ([v^*v])$,
$g=\varphi ([1-v^*v])$ and $h=\varphi ([1-vv^*])$, which are functions
in $\W$. Note that $f+g=f+h=d$, and that both $g_{|\de S_u}$ and
$h_{|\de S_u}$ are finite functions, since $g,h\in
I_{fin}:=\varphi(V(I_{fin}(A)))$. If $g(s)=h(s)$, then $g=h$ and so
$1-v^*v\sim 1-vv^*$ in $I$, whence $\del [w]=0$. Suppose that
$g(s)<h(s)$. Then, by Lemma~\ref{2.3} $(a)$ there exist $e\in \LA
(S_u)^{++}$ and $k\in \mathbb{N}$ such that $g+k=h+e$ and
$e+d=k+d$. Thus there exist nonzero projections $p\in \Moo (L(A))$ and
$q\in \Moo (I)$ such that $\varphi ([p])=k$ and $\varphi ([q])=e$, and
$(1-v^*v)\oplus p\sim (1-vv^*)\oplus q$. Hence $\pi (1-v^*v)\sim \pi
(1-vv^*)\oplus \pi (q)$ in $I/L(A)$, and it follows that $\del ([w])
=[\pi (1-v^*v)]-[\pi (1-vv^*)]\geq 0$. A similar argument, if
$h(s)<g(s)$, shows that $\del ([w])\leq 0$.

We conclude from Proposition~\ref{2.2} that $w$ can always be lifted
either to an isometry or to a co-isometry.\qed

If $X$ is a subset of a convex set $K$, we denote by ${\rm conv}(X)$
the convex hull of $X$, that is, the smallest convex subset of $K$
that contains $X$.
\begin{lem}
\label{2.4a}
Let $K$ be a Choquet simplex, and let $F$ be a closed face of $K$. Let
$F'$ be the complementary face of $F$. Then $\ol{F'}=\ol{{\rm
    conv}(\de K\setminus\de F)}$.
\end{lem}
{\bf Proof.}
Suppose first that $F'$ is a closed face of $K$. Then, since $K$ is
compact, it follows that $F'$ is a compact convex subset of $K$. By
Krein-Mil'man's Theorem, we have that $F'=\ol{{\rm conv}(\de
  F')}$. Note now that $\de F'=F'\cap \de K$, because $F'$ is a face
of $K$, and by definition of complementary face we conclude that $\de
F'=\de K\setminus \de F$.

If $F'$ is not closed, set $X:=\ol{{\rm conv} (\de K\setminus\de
  F)}$. By construction, $X$ is a compact convex subset of $K$, and it
contains all the extreme points of $K$, except maybe those from
$F$. Notice now that the convex hull of $X$ and $F$ is closed (by
\cite[Proposition 5.2]{poag}), and it contains all the extreme points
of $K$. Another application of Krein-Mil'man shows that $K={\rm
  conv}(X\cup F)$. Let $a\in F'$. Then there exist $\alp\in [0,1]$,
$f\in F$ and $x\in X$ such that $a=\alp f+(1-\alp)x$. Since $f\notin
F'$, and taking into account that $F'$ is a face, we have that
$\alp=0$, and thus $a=x$. We therefore conclude that $F'\subseteq
X$. It is clear, on the other hand, that $X\subseteq \ol{F'}$, whence
it follows that $\ol{F'}=X$, as desired.\qed

\begin{prop}
\label{2.5}
Let $A$ be a (nonunital) simple separable \C\ with real rank zero and
stable rank one. Suppose that $A$ is nonelementary and that $V(A)$ is
strictly unperforated. Let $p\in A$ be a nonzero projection, and
assume that $A$ has exactly one infinite extremal quasitrace in $Q_p$,
which we denote by $\tau$. If $RR(\mul)=0$, then all proper isometries
of $\mul/I_{fin}(A)$ can be lifted to proper isometries of $\mul/L(A)$
if and only if the complementary face of $\{\tau\}$ in $Q_p$ is
closed.
\end{prop}
{\bf Proof.}
Let $u=[p]\in V(A)$, and set $d=\sup\phi_u(D(A))$. By
Theorem~\ref{tm}, there is a monoid isomorphism $\varphi: V(\mul)\fl
V(A)\sqcup \W$ such that $\varphi([1_{\mul}])=d$. As in the previous
result, the quasitrace $\tau$ corresponds, through the affine
homeomorphism $\alp:Q_p\fl S_u$, to a unique point $s\in \de S_u$ such
that $d(s)=\infty$.

First, suppose that $\{s\}'$, the complementary face of $\{s\}$ in the
simplex $S_u$, is not closed and that all proper isometries of $\mul
/I_{fin}(A)$ can be lifted to proper isometries of $\mul /L(A)$. Since
all projections in $\mul /I_{fin}$ are infinite (see
Theorem~\ref{tm3}), there exists a proper isometry $w\in \mul
/I_{fin}(A)$. By hypothesis, $w$ can be lifted to a proper isometry
$v\in \mul /L(A)$. On the other hand, there exists a partial isometry
$z\in \mul $ such that $\pi (z)=v$, where $\pi:\mul \fl \mul /L(A)$ is
the natural map (again by \cite[Proof of Lemma 2.6]{elan} or
\cite[Lemma 2.1]{pext}). Therefore $1-z^*z\in L(A)$ and $1-zz^*\notin
I_{fin}(A)$. As in Proposition~\ref{2.4}, we have that $z\notin A$,
and moreover in this case $1-zz^*\notin A$. We may also assume that
$1-z^*z\notin A$, since the proof would be similar otherwise. Let
$f=\varphi ([z^*z])$, $g=\varphi ([1-z^*z])$ and $h=\varphi
([1-zz^*])$. Then $f+g=f+h=d$, and $h(s)=\infty$. Note also that $g\in
\varphi(V(L(A)))=V(A)\sqcup \Aff (S_u)^{++}$, and thus is
continuous. Since $d_{|\de S_u}$ is infinite exactly at $s$ we have
that $f(s)=\infty$. Therefore $g_{|\de S_u\setminus\{s\}}=h_{|\de
  S_u\setminus\{s\}}$. We also have that $F:=\ol{\{s\}'}$ equals the
closed convex hull of $\de S_u\setminus\{s\}$, by Lemma~\ref{2.4a},
and that $F$ is a compact convex subset of $S_u$ (see
\cite[Proposition 5.1]{poag}).

If $s\notin F$, then $\de F=\de S_u\setminus\{s\}$, and since $g$ and
$h$ are affine and lower semicontinuous we get that $g_{|F}=h_{|F}$
(by Lemma~\ref{pacha}). Since $\{s\}'$ is not closed, there exists
$x\in F\setminus\{s\}'$. Thus $g(x)=h(x)$. But $x\notin \{s\}'$, so
that there exist $\alp \in (0,1]$ and $t\in \{s\}'$ such that $x=\alp
s+(1-\alp )t$, and this implies $h(x)=\infty$, a contradiction since
$g$ is continuous.

Hence $s\in F$, and so $s=\lim\limits_n y_m$, where $y_m$ belong (for
all $m$) to the convex hull of $\de S_u\setminus\{s\}$. Thus
$g(y_m)=h(y_m)$ for all $m$, and it follows that $g(s)=\lim\limits_m
g(y_m)=\lim\limits_m h(y_m)\geq h(s)$, a contradiction.

Conversely, assume that $\{s\}'$ is a closed face, and let $w\in \mul
/I_{fin}(A)$ be a proper isometry. Again, there exists a partial
isometry $z\in \mul \setminus A$ such that $\phi (z)=w$, where
$\phi:\mul \fl \mul /I_{fin}(A)$ is the natural map. Thus $1-z^*z \in
I_{fin}(A)$ and $1-zz^*\notin I_{fin}(A)$. If $1-z^*z\in A$, then
$1=\pi (z)^*\pi (z)$ and $1\ne \pi (z)\pi (z)^*$ (where $\pi$ is the
quotient map modulo $L(A)$). We also have that $w=\phi(z)=\ol{\pi
  (z)}$, where the latter denotes the class of $\pi(z)$ modulo
$I_{fin}(A)/L(A)$. We then conclude that $w$ can be lifted to a
(proper) isometry of $\mul/L(A)$, by Lemma~\ref{2.1}. Therefore we may
assume that $1-z^*z\notin A$. Let $f=\varphi ([z^*z])$, $g=\varphi
([1-z^*z])$ and $h=\varphi ([1-zz^*])$. Then $f+g=f+h=d$ and also
$h(s)=\infty $. By Lemma~\ref{2.3} $(b)$, there exist $e\in \LA
(S_u)^{++}$ and $k\in \mathbb{N}$ such that $g+e=h+k$ and
$e+d=k+d$. Hence, there exist projections $p\in \Moo (L(A))$ and $q\in
\Moo (\mul )$ such that $k=\varphi ([p])$ and $e=\varphi ([q])$, and
$(1-z^*z)\oplus q\sim (1-zz^*)\oplus p$. Therefore $1-\pi (z)^*\pi
(z)\lesssim 1-\pi (z)\pi (z)^*$, and as before $w=\phi(z)=\ol{\pi
  (z)}$; it follows then from Lemma~\ref{2.1} that $w$ can be lifted
to a (proper) isometry.\qed

We turn our attention now to the extreme points of $\mul /L(A)$,
which, in our setting, are only isometries and co-isometries.
\begin{lem}
\label{2.6}
Let $A$ be a (nonunital) simple separable \C\ with real rank zero and
stable rank one. Suppose that $A$ is nonelementary and that $V(A)$ is
strictly unperforated. Let  $p\in A$ be a nonzero projection, and
assume that $A$ has exactly one infinite extremal quasitrace in
$Q_p$. If $L(A)$ has real rank zero, then $\goe (\mul/L(A))$ consists
only of isometries and co-isometries.
\end{lem}
{\bf Proof.}
Let $u=[p]\in V(A)$ and $d=\sup\phi_u(D(A))$. As before, we have a
monoid isomorphism $\varphi:V(\mul )\fl V(A)\sqcup \W$ such that
$\varphi([1_{\mul}])=d$, and there is a unique state $s\in \de S_u$ at
which $d_{|\de S_u}$ is infinite.

Let $v\in \goe (\mul /L(A))$. Suppose that $v$ is neither an isometry
nor a co-isometry. Choose a partial isometry $z\in \mul \setminus A$
such that $\pi (z)=v$, where $\pi:\mul \fl \mul /L(A)$ is the natural
map. Then $(1-z^*z)\mul (1-zz^*)\subseteq L(A)$ and
$1-z^*z,1-zz^*\notin L(A)$. Let $I$ and $J$ denote the closed ideals
of $\mul$ generated by $1-z^*z$ and $1-zz^*$ respectively. Then
$IJ\subseteq L(A)$. On the other hand, we have that $A\subsetneqq I,J$
since $v$ is neither an isometry nor a co-isometry. Therefore
$L(A)\subseteq I,J$ and we conclude that $L(A)=IJ=I\cap J$.

Clearly $1-z^*z,1-zz^*\notin A$. Let $f=\varphi ([z^*z])$, $g=\varphi
([1-z^*z])$ and $h=\varphi ([1-zz^*])$. Then again $f+g=f+h=d$, whence
$g=h$ on $\de S_u\setminus\{s\}$. Let
$$I_g=V(A)\sqcup \{f_1\in\W\mid f_1+f_2=ng \mbox{ for some }f_2\in\W
\mbox{ and some }n\in\mathbb{N}\}$$
be the order-ideal of $V(A)\sqcup\W$ generated by $g$. Similarly,
$I_h$ denotes the order-ideal generated by $h$. Note that
$\varphi(V(I))=I_g$ and $\varphi(V(J))=I_h$. Since $I\cap J=L(A)$, we
have that
$$I_g\cap I_h=\varphi(V(I)\cap V(J))=\varphi(V(L(A)))=V(A)\sqcup \Aff
(S_u)^{++}.$$
Note that $g,h\notin \Aff (S_u)^{++}$, since $v$ is neither an
isometry, nor a co-isometry. Therefore $g\neq h$. (Otherwise
$I_g=I_g\cap I_h=V(A)\sqcup \Aff (S_u)^{++}$, and so $g$ would be
continuous, a contradiction.) In particular, $g$ and $h$ are not
simultaneously infinite. Suppose that $g(s)=\infty $ and $h(s) <\infty
$. Then $g+h=2g$, whence $I_h\subseteq I_g$, so that $h$ would be
continuous, a contradiction. The argument is similar if $g(s)<\infty$
and $h(s)=\infty$.

This implies that $g(s),h(s)<\infty$. Suppose that $h(s)<g(s)$. Let
$n\in \mathbb{N}$ be such that $n\geq 2$ and $nh(s)>g(s)$.  By a
similar argument to the one used in \cite[Proposition 4.13]{id}, there
exists $f'\in \LA (S_u)^{++}$ such that $f'(s)=nh(s)-g(s)$ and
${f'}_{|\{ s\}'}=(n-1)g_{|\{s\}'}$. Then $f'+g=nh$, and therefore
$I_g\subseteq I_h$, a contradiction since $g$ is not
continuous. Again, the argument is similar if $g(s)<h(s)$.

Hence we conclude that any $v\in \goe (\mul /L(A))$ is necessarily an
isometry or a co-isometry.\qed

\begin{rema}
{\rm The previous result would follow immediately if $\mul/L(A)$ were
  a prime ring. We remark that this is not true in general.}
\end{rema}
{\bf Proof.}
Let $A$ be a (nonunital) simple separable $AF$ algebra, and let $p\in
A$ be a nonzero projection such that if $u=[p]\in V(A)$, then $\de
S_u\cong [-1,1]$. Moreover, we take $A$ such that its scale $d$ equals
$2$ on the interval $[-1,0)$, and such that $d(x)=1/(1-x)$, for $x\in
[0,1]$. (The existence of this example follows after \cite[Example
4.7]{id}; see also \cite[Example 7.3]{gok}.) Then $A$ has only one
infinite extremal (quasi)trace. Define $f\in\W$ as $1$ on $[-1,1)$,
and set $f(1)=1/2$. Define also $g\in\W$ as $1$ on $[-1,0)$ and $1/2$
on $[0,1]$. Then, if we denote by $I_f$ and $I_g$ the respective
order-ideals generated by $f$ and $g$, it is clear that $I_f\cap
I_g=V(A)\sqcup C[-1,1]^{++}$, whereas $I_f\ne V(A)\sqcup C[-1,1]^{++}$
and also $I_g\ne V(A)\sqcup C[-1,1]^{++}$. Since $\de S_u$ is compact,
it follows from Theorem~\ref{tm} and \cite[Theorem 2.3]{zhriesz} (see
also \cite[Section 6]{id}) that $\mul/L(A)$ is not a prime ring.\qed

Recall that a closed ideal $I$ of a \C\ $B$ is {\bf stably cofinite}
provided that the quotient $B/I$ is stably finite. This is equivalent
to saying that the monoid $V(B/I)$ is stably finite, that is, the
relation $x+y=y$ implies $x=0$. When studying the stably finite
quotients of $\mul$, a closed ideal called the bounded ideal is of
some significance. In our setting, it can be described as the unique
closed ideal $I_b(A)$ of $\mul$ such that
$\varphi(V(I_b(A)))=V(A)\sqcup\{f\in\W\mid f \mbox{ is bounded}\}$,
where $\varphi$ is the isomorphism established on Theorem~\ref{tm}. If
$A$ has a finite number of infinite extremal quasitraces, it is
possible to characterize when this ideal is stably cofinite, as
follows.
\begin{prop}
\label{sc}
Let $A$ be a (nonunital) simple separable \C\ with real rank zero and
stable rank one. Suppose that $A$ is nonelementary, that $V(A)$ is
strictly unperforated and that $\mul$ has real rank zero. Let $p\in A$
be a nonzero projection, and denote by $Q_{\infty}$ the set of
infinite extremal quasitraces in $Q_p$. If $Q_{\infty}$ is finite,
then $I_b(A)$ is stably cofinite if and only if for any nonempty
subset $X\subseteq Q_{\infty}$, the complementary face of ${\rm
  conv}(X)$ (in $Q_p$) is not closed.
\end{prop}
{\bf Proof.}
Let $u=[p]\in V(A)$, and let $d=\sup\phi_u(D(A))$. Suppose that
$Q_{\infty}$ has cardinality $n$. Since the natural map $\alp:Q_p\fl
S_u$ is an affine homeomorphism, the fact that $Q_{\infty}$ is finite
means exactly that the set $\Gamma_d:=\{s\in\de S_u\mid d(s)=\infty\}$
is finite, and with cardinality $n$. Set
$\Gamma_d=\{s_1,\ldots,s_n\}$.

Suppose that $I_b(A)$ is stably cofinite. Then the order-ideal
$I_b=\varphi(V(I_b(A)))$ of $V(A)\sqcup\W$ is stably cofinite, where
$\varphi$ is the isomorphism in \ref{tm}. Suppose that there exists a
nonempty subset $X\subseteq\Gamma_d$ such that $F:=({\rm conv}(X))'$,
the complementary face of ${\rm conv}(X)$, is closed. Then we have
$\Aff (S_u)\cong \Aff ({\rm conv}(X))\times \Aff(F)$ (\cite[Corollary
11.27]{poag}). We may define functions $f_n\in\Aff (S_u)^{++}$ by
${f_n}_{|{\rm conv}(X)}=n$ and ${f_n}_{|F}=1$. Then, let
$f=\sup\limits_nf_n$. We have that $f\in {\rm LAff}(S_u)^{++}$ and
that $f_{|X}=\infty$, while $f_{|F}=1$. Thus $f+d=1+d$, so that
$f\in\W$, and in the quotient $(V(A)\sqcup\W)/I_b$ we have
$[f]+[d]=[d]$ with $[f]\ne 0$, which contradicts the stable finiteness
of $I_b$.

Conversely, suppose that for any nonempty subset $X\subseteq\Gamma_d$,
the complementary face of ${\rm conv}(X)$ is not closed. Suppose that
there exist bounded functions $l_1,l_2\in\W$, a number
$n\in\mathbb{N}$ and a function $g\in\W$ such that
$nd+l_1=nd+g+l_2$. Then $g$ is bounded on $\de S_u\setminus\Gamma_d$,
say $g_{|\de S_u\setminus\Gamma_d}\leq M$. Let $F=({\rm
  conv}(\Gamma_d))'$. Since $F$ is not closed, there exists $x\in
\ol{F}\setminus F$. By Lemma~\ref{2.4a}, there is a sequence $\{t_k\}$
in ${\rm conv}(\de S_u\setminus\Gamma_d)$ that converges to $x$. Since
$g$ is affine, we have that $g(t_k)\leq M$ for all $k$, whence
$g(x)\leq M$, due to lower semicontinuity. On the other hand, taking
into account that $S_u$ is the direct convex sum of ${\rm
  conv}(\Gamma_d)$ and $F$, there exist positive numbers $\alp_i$ for
$i=1,\ldots,n$ (not all zero), a number $\beta\geq 0$ and an element
$t\in F$ such that $\sum\limits_{i=1}^n \alp_i+\beta=1$, and such that
$x=\sum\limits_{i=1}^n\alp_is_i+\beta t$. Therefore there exists
$1\leq i\leq n$ such that $g(s_i)<\infty$. Without loss of generality,
we may assume that $g(s_1)<\infty$. Now, a recursive argument shows
that $g(s_2),\ldots,g(s_n)<\infty$, whence we conclude that $g$ is
bounded, and thus $I_b$ is stably cofinite.\qed

We are now in position to characterize the extremal richness of the
corona algebra for a simple \C\ with exactly one infinite extremal
quasitrace. This situation was considered in \cite[Proposition
4.18]{lo}, for separable, nonunital and simple $AF$ algebras with
finitely many extremal traces. In that case the corona algebra is
always extremally rich. As we will see, the general situation is quite 
different.
\begin{theor}
\label{2.7}
Let $A$ be a (nonunital) simple separable \C\ with real rank zero and
stable rank one. Suppose that $A$ is nonelementary and that $V(A)$ is
strictly unperforated. Let $p\in A$ be a nonzero projection, and
assume that $A$ has exactly one infinite extremal quasitrace in $Q_p$,
which we denote by $\tau$. If $RR(\mul)=0$, then the following
conditions are equivalent:
\casos
\cas
$\mul/A$ is extremally rich;
\cas
The complementary face of $\{\tau\}$ (in $Q_p$) is closed;
\cas
The ideal $I_b(A)$ is not stably cofinite.
\fcasos
\end{theor}
{\bf Proof.}
Note first that $L(A)/A$ is a closed essential ideal of $\mul /A$,
which is simple and purely infinite. Then by Theorem~\ref{1.3a}, $\mul
/A$ is extremally rich if and only if $\mul /L(A)$ is extremally rich
and $\goe (\mul /L(A))$ consists only of isometries and
co-isometries. The latter condition is automatic from
Lemma~\ref{2.6}. Thus $\mul /A$ is extremally rich if and only if
$\mul /L(A)$ is extremally rich.

Since $I_{fin}(A)/L(A)$ has stable rank one (see \cite[Proposition
6.1]{id}), it follows from Corollary~\ref{1.1a} that $\mul /L(A)$ is
extremally rich if and only if $\mul /I_{fin}(A)$ is extremally rich
and the extreme partial isometries of $\mul /I_{fin}(A)$ can be lifted
to those of $\mul /L(A)$. Note that $\mul /I_{fin}(A)$ is purely
infinite and simple (by \cite[Theorem 6.3, Proposition 6.5]{id}),
hence extremally rich by Remark~\ref{rm} $(b)$. Therefore $\mul /A$ is
extremally rich if and only if the extreme partial isometries of $\mul
/I_{fin}(A)$ can be lifted to those of $\mul /L(A)$, and this last
condition holds if and only if the complementary face of $\{\tau\}$ in
$Q_p$ is a closed face, by Proposition~\ref{2.4} and
Proposition~\ref{2.5}. This proves $(a)\sns (b)$.

The equivalence between $(b)$ and $(c)$ follows directly from
Proposition~\ref{sc}.\qed

\begin{corol}
\label{2.8}
Let $A$ be a (nonunital) simple separable \C\ with real rank zero and
stable rank one. Assume that $A$ is nonelementary and that $V(A)$ is
strictly unperforated. Let $p\in A$ be a nonzero projection, and
suppose that $A$ has exactly one infinite extremal quasitrace $\tau$
in $Q_p$. If $RR(\mul)=0$ and $\de Q_p$ is a compact space, then
$\mul/A$ is extremally rich if and only if $\tau$ is isolated in $\de
Q_p$.
\end{corol}
{\bf Proof.}
By Theorem~\ref{2.7} we have to prove that $\{\tau\}$ is isolated in
$\de Q_p$ if and only if the complementary face of $\{\tau\}$ in $Q_p$
is closed.

Let $u=[p]\in V(A)$ and $d=\sup\phi_u(D(A))$. Let $s\in \de S_u$ be
the unique extreme point at which $d_{|\de S_u}$ is infinite. This
point exists because there is an affine homeomorphism $\alp:Q_p\fl
S_u$ (and in fact $s:=\alp(\tau)$). We have to prove then that $\{s\}$
is isolated in $\de S_u$ if and only if $\{s\}'$ is a closed
face. Note that $\{s\}$ is isolated if and only if $\de
S_u\setminus\{s\}$ is closed.

By \cite[Corollary 11.20]{poag} (and since $\de S_u$ is compact),
there is an affine homeo\-mor\-phism $\beta$ between $S_u$ and
$M_1^+(\de S_u)$ given by $\beta(t)=\eps_t$, for $t\in \de S_u$, and
where $\eps_t$ is the point mass measure at $t$.

Suppose that $\{\eps_s\}'$ is a closed face in $M_1^+(\de S_u)$. It
follows from \cite[Proposition 5.25]{poag} that there exists a closed
subset $X\subseteq\de S_u$ such that $\{\eps_s\}'=\sigma(X):=\{\mu\in
M_1^+(\de S_u)\mid \mu(X)=1\}$. In fact, $X=\{t\in \de S_u\mid
\eps_t\in \{\eps_s\}'\}$, which clearly equals $\de
S_u\setminus\{s\}$, and therefore $\de S_u\setminus\{s\}$ is closed.

Conversely, if $X=\de S_u\setminus\{s\}$ is closed, letting $F=\sigma
(X)$, we have as before that $F$ is a closed face, which is easily
shown to coincide with $\{\eps_s\}'$. Thus we have proved that
$\{s\}'$ is a closed face if and only if $\de S_u\setminus\{s\}$ is
closed, and the result follows.\qed

The next instance of interest is to determine when $\mul$ has extremal
richness. In \cite[Theorem 3.1]{lo} it is proved that if $A$ is a
separable unital simple \C\ with real rank zero and with a finite
trace, then $\mathcal{M}(A\otimes \mathbb{K})$ is not extremally
rich. In the same paper it is noted that $\mul$ is not extremally rich
when $A$ is a simple, separable $AF$ algebra without unit and
nonelementary. Part of the argument in the following is based on
\cite[Proof of Theorem 2.3]{lo}, which we include for completeness.
\begin{prop}
Let $A$ be a (nonunital) simple separable \C\ with real rank zero,
stable rank one and with $V(A)$ strictly unperforated. Assume also
that $A$ is nonelementary. If $RR(\mul)=0$, then $\mul$ is not
extremally rich.
\end{prop}
{\bf Proof.}
Suppose first that the scale of $A$ is not identically infinite. Then,
by \cite[Corollary 7.8]{id} $\mul$ is stably finite. Notice also that
$\mul$ is a prime 
ring and thus $\goe (\mul)$ consists of isometries and
co-isometries. We conclude that $\goe (\mul)=\mathcal{U}
(\mul)$. According to Remark~\ref{rm} $(a)$ we have that $\mul$ is
extremally rich if and only if ${\rm sr}(\mul)=1$. Since $V(\mul)$ is
not cancellative, we conclude that ${\rm sr}(\mul)\ne 1$, whence
$\mul$ is not extremally rich.

Suppose now that the scale of $A$ is identically infinite. If $A$ has
at least two extremal quasitraces (which will be infinite), then
$\mul$ is not extremally rich, by Theorem~\ref{1.4} $(b)$. We are
therefore left to consider the case when $A$ has only one extremal
quasitrace, which is infinite. This means that if $\{e_n\}$ is an
(increasing) approximate unit for $A$ consisting of projections, then
$\sup\tau(e_n)=\infty$. Let $\pi:\mul\fl\mul/A$ be the natural
quotient map, and let $q\in L(A)\setminus A$ be a projection. Then
$\pi (q)$ is an infinite projection, and therefore there exists $v\in
L(A)/A$ such that $v^*v=\pi (q)$ while $vv^*<\pi (q)$.

Let $w=v+1-\pi (q)$, a proper isometry of $\mul/A$ which lifts (since
$RR(A)=0$ and by the proof of \cite[Lemma 2.6]{elan}) to a partial
isometry $u\in\mul$. There exists (by an argument used in
\cite[Theorem 2.2]{lo}) a projection $p\in A$ such that $u^*u=1-p$. It
follows that $1-uu^*\in L(A)$ and if $\{f_n\}$ is an approximate unit
of projections for the $C^*$-subalgebra $(1-uu^*)A(1-uu^*)$ of $A$, we
have that $\sup\tau (f_n)<\infty$.

Let $p'\in A$ be a projection such that $\tau
(p')>\sup\tau(f_n)$. Since $1-p,1-p'\notin L(A)$, and since $A$ has
only one extremal quasitrace (that is infinite) we get that $1-p\sim
1-p'$. Thus there is $r\in\mul$ such that $1-p'=r^*r$ while
$1-p=rr^*$. Let $t=ur$. Then $t^*t=1-p'$ and $tt^*=uu^*$, so that $\pi
(t)$ is an isometry in $\mul/A$.

Denote by $\phi:\mul\fl\mul/L(A)$ the natural quotient map, and note
that in fact $\phi(t)$ is a unitary in $\mul/L(A)$. We claim that
$\pi(t)$ cannot be lifted to any isometry of $\mul$. If there exists
an isometry $s\in\mul$ such that $\pi(s)=\pi(t)$, then
$\phi(s)=\phi(t)$, and therefore
$$\del([\phi(t)])=[1-s^*s]-[1-ss^*]=[1-t^*t]-[1-tt^*],$$
where $\del:K_1(\mul/L(A))\fl K_0(L(A))$ is the index map. It follows
that $[1-tt^*]=[1-t^*t]+[1-ss^*]$. Since $\mul$ has real rank zero, we
get from Lemma~\ref{g} a projection $f\in\Moo (L(A))$ such that
$$(1-tt^*)\oplus f\sim (1-t^*t)\oplus (1-ss^*)\oplus f=p'\oplus
(1-ss^*)\oplus f.$$
By simplicity, we assume that $f\in L(A)$. Let $\{t_n\}$
(resp. $\{g_n\}$) be an approximate unit of projections for
$(1-ss^*)A(1-ss^*)$ (resp. for $fAf$). Then by \cite[Proposition
1.7]{gok}, for each $n\in \mathbb{N}$, there exists $m\in\mathbb{N}$
such that $p'\oplus t_n\oplus g_n\lesssim f_m\oplus g_m$. Therefore,
if $n\in\mathbb{N}$, there is $m\in\mathbb{N}$ such that $\tau
(p')+\tau (t_n)+\tau(g_n)\leq \tau(f_m)+\tau(g_m)$. It follows that
$\sup\tau (t_n)\leq 0$, and this implies that $1-tt^*=0$. Thus
$1-uu^*=0$ and hence $ww^*=1$. This contradicts the fact that $w$ is a
proper isometry, and therefore the claim is established.

It is easy to see that $\pi(t)$ cannot be lifted to a co-isometry,
either. We then conclude from Theorem~\ref{1.1} that $\mul$ is not
extremally rich.\qed

\section*{Acknowledgments}

The author wishes to thank Pere Ara for many helpful comments and
for allowing the inclusion of Lemma~\ref{g}.

This research has been partially supported by MEC-DGICYT grant
no.PB95-0626, and by the Comissionat per Universitats i Recerca de la
Generalitat de Catalunya. This is part of the author's Ph.D. Thesis,
written under the supervision of Professor Pere Ara.

\end{document}